\documentclass[12pt,letterpaper,fleqn]{article}
\usepackage{amsfonts}
\usepackage{graphicx}
\usepackage{amsmath}
\newcommand{\D}{\displaystyle}
\newcommand{\FS}{\scriptsize}
\newcommand{\dsum}{\displaystyle \sum}
\newcommand{\Dfrac}{\displaystyle \frac}
\newcommand{\dint}{\displaystyle \int}

\begin{document}

\title{\large Computational Techniques for Simulating Natural
Convection  in Three-Dimensional Enclosures  with Tetrahedral
Finite Elements} \author{\large K. O. Ladipo \and R.
Glowinski\and T.W. Pan}
\date{ }
\maketitle

\begin{abstract}
\setlength{\parindent}{0pt}
This article  discusses computational techniques  for
simulating natural convection in  three-dimensional domains using
finite element methods with tetrahedral elements. These
techniques form a new numerical procedure for this kind of
problems.  In this procedure, the treatment of advection by a
\emph{wave equation approach} is extended to three-dimensional
unstructured meshes with tetrahedra.

Numerical results of  natural convection of an
 incompressible  Newtonian fluid in a cubical enclosure at
Rayleigh numbers in the range $10^3$ to $10 ^6$  are obtained
and they are in good agreement with those in literature obtained
by other methods.
 
\textbf{key words} Three-dimensional domains, Finite element
methods, Tetrahedral elements, Natural convection, Incompressible
fluids.
\end{abstract}
 
\pagenumbering{arabic}

\section {\normalsize Introduction}

Simulation of natural convection flows in three-dimensional
geometries has been an area of active research in recent years.
In the  past decade, most researchers  who  performed
calculations  in three-dimensional  geometries   were hindered
from applying sufficient  resolutions, by limitations on computer
storage. For example, Mallinson and De Vahl   Davis in
\cite[1977]{bib:mallinsondv77}  used a very coarse mesh with up
to $15^3$ nodes, Pepper D.W.  in \cite[1987]{bib:pepper87}
applied only $33\times17\times9$ nodes. More recently, with  the
availability of  more powerful computers,    researchers are
now able to perform  calculations on meshes with better
resolutions.   Le Peutrec and Lauriat \cite[1990]{bib:peutrec90}
used  mesh with up to $41^3$  nodes,   Fusegi \emph{et al.\@} in
\cite[1991]{bib:fusegi91}  used meshes with up to $62^3$  nodes.
Janssen \emph{et al.\@}in  \cite[1993]{bib:janshen93} reported
results with $120^3$ nodes mesh,  but they generated the results
by  symmetry,  after performing  actual simulation  with  only
one-fourths of this number.
 
Natural convection is governed by a coupled system of
Navier-Stokes equations and energy equations.

The objective of this paper is to present a finite element method
 for simulating  natural convection  of an incompressible fluid
in three-dimensional geometries using tetrahedral elements with
unstructured mesh. An operator-splitting scheme of
Marchuk-Yanenko is applied to split the coupled system into three
sub-problems namely, the pressure, transport and diffusion
sub-problems. This decouples the  difficulties usually associated
 with non-linearity and incompressibility constraint. The pressure and diffusion sub-problems are time
discretized  by backward-Euler-type method. The non-linear
advection is treated by a \emph{wave equation approach}. Space
discretization is achieved by a finite element  method where
pressure,  velocity and temperature are approximated by
continuous piecewise-linear polynomials  on meshes consisting of
4-node tetrahedral elements.  The mesh for velocity and
temperature is twice finer  than the pressure mesh so that the
\emph{inf-sup} condition  is satisfied. A systematic method of
constructing these velocity-pressure meshes, such that  each
pressure tetrahedral element is a \emph{macro-element} consisting
of eight sub-tetrahedra for velocity, is discussed in this
article. We extend the  two-dimensional method  for constructing
a pressure macro-element, by connecting edge mid-points, to three
dimensional meshes. The numerical procedure presented in this
article also extends the treatment of advection by a \emph{wave
equation approach} in  \cite{bib:Mythesis}, \cite{bib:DrPanDrGlo}
to  three-dimensions while combining other different numerical
techniques  and thus forming  a  new, efficient, solution
procedure suitable for simulating motion of an
incompressible fluid in three-dimensional geometries  with
unstructured meshes.

Results obtained for the  numerical example of natural convection
 of air in a cubical box, illustrate  the accuracy and
reliability of this new procedure. The three-dimensional results
also validate the  usual  assumptions in two-dimensional
simulations and elucidates three-dimensional effects on this flow
phenomenon.

\section{\normalsize Governing Equations  for Natural Convection}

We consider natural convection of an incompressible viscous
Newtonian  fluid  enclosed in a three-dimensional rectangular
domain, $\Omega \subset \mathbf{R}^3$, with  boundary denoted by
$\Gamma$. The geometry and  coordinate system for  the  enclosure
 are shown in Figure~\ref{fig:geom3D}.  The natural convection is
induced by the non-zero temperature gradient between the two
vertical surfaces,  $\Gamma_l$ (at $x = 0$) and $\Gamma_r$ (at $x
= L_x$). The remaining four surfaces, $(\Gamma\setminus
\Gamma_l\cup \Gamma_r)$ are assumed to be perfectly thermally
insulated.

\setlength{\unitlength}{0.8cm}
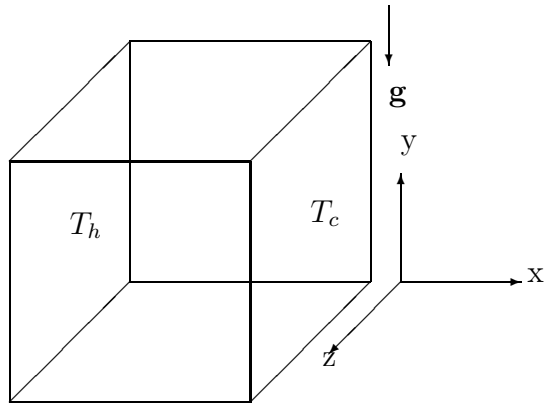
\begin{figure}[ht] \centering
\begin{picture}(7,7)(0,0)
\put(0,0){\line(1,0){4}}
\put(4,0){\line(0,1){4}}
\put(4,4){\line(-1,0){4}}
\put(0,4){\line(0,-1){4}}
\put(2,2){\line(1,0){4}}
\put(6,2){\line(0,1){4}}
\put(6,6){\line(-1,0){4}}
\put(2,6){\line(0,-1){4}}
\put(0,0){\line(1,1){2}}
\put(4,0){\line(1,1){2}}
\put(4,4){\line(1,1){2}}
\put(0,4){\line(1,1){2}}
\put(6.3,6.6){\vector(0,-1){1}}
\put(6.3,5){\textbf{g}}
\put(6.5,2){\vector(0,1){1.8}}
\put(6.5,2){\vector(1,0){2}}
\put(8.6,2){{x}}
\put(6.5,4.2){{y}}
\put(6.5,2){\vector(-1,-1){1.2}}
\put(5.2,0.6){{z}}
\put(1.0,2.8){$T_h$}
\put(5.0,3.0){$T_c$}
\end{picture}
\caption[Geometry and coordinate system for cubical enclosures]{Geometry and coordinate system for
cubical enclosures,\\ \mbox{}\hspace{0.8in}( $0 \leq x \leq L_x \ ; \  0 \leq y \leq L_y \ ; \  0 \leq z \leq L_z$).}\label{fig:geom3D}
\end{figure}

With the Boussinesq approximation, the  vector form of the dimensionless
governing equations in a finite time interval $(0,\, t^N)$ are:

\begin{eqnarray}
\frac{\partial \mathbf{u}}{\partial t} +  (\mathbf{u} \cdot
\nabla) \mathbf{u} - Pr \Delta \mathbf{u} +  \nabla p\ =\  Ra Pr \theta \hat{\jmath}\ \
\mbox{in}\  \Omega \times (0,t^N),& &  \label{eq:NDNS}\\
\frac{\partial \theta}{\partial t}\ +\   (\mathbf{u} \cdot \nabla)
\theta\ -\  \Delta \theta\  =\  0\hspace{1in}  \mbox{in}\ \ \Omega \times (0 ,
t^N),& & \label{eq:NDheat}\\
\nabla \cdot \mathbf{u}\  =\  0\hspace{1in}  \mbox{in}\ \ \Omega \times (0 ,
t^N),& &   \label{eq:NDdiv}
\end{eqnarray}
where ,
\begin{gather}Ra = \frac{g \beta (T_h - T_c) L_x^3}{\nu
\alpha},\mbox{  the Rayleigh number},\\
Pr = \frac{\nu}{\alpha},\mbox{  the Prandtl number},
\end{gather}
\(\mathbf{u}( x,\,y,\,z,\,t)\ =\ \left\{ u_i\right\}^3_{i = 1} =  (u_x,\, u_y,\, u_z)\)  is  the flow velocity,

$t$ is elapsed time,
\( (\mathbf{u} \cdot \nabla )\mathbf{u}\) is a symbolic notation
for the non-linear vector term
 $ \left\{ \dsum^3_{j = 1} u_j \frac{\partial u_i}{\partial x_j} \right\}^3_{i = 1}$,

$\left\{ x_i\right\}^3_{i = 1}  = ( x,\,y,\,z )$,

$\theta( x,\,y,\,z,\,t)$ is the temperature,

$p( x,\,y,\,z,\,t)$ is the pressure ,

$\beta$ is the coefficient of thermal expansion of the fluid,

$\alpha$ is the coefficient of thermal diffusivity of the fluid,

$\nu$ is the kinematic viscosity coefficient of the fluid,

$\mathbf{g}$ is the gravitational acceleration.

This set of dimensionless equations is subject to the following \emph{initial and boundary conditions}:

\begin{equation}\left\{
\begin{array}{l}                                        1
\mathbf{u}(x,\,y,\,z,\, 0)\ =\ \mathbf{0}\hspace{0.2in}\mbox{in}\ \
{\Omega},\\ \theta( x,\,y,\,z,\, 0)\ =\  0
 \hspace{0.2in}\mbox{in}\ \ {\Omega},\\ \mathbf{u}(x,\,y,\,z,\,
t)\ =\ \mathbf{0}\hspace{0.2in}\mbox{on}\ \ \Gamma\times (0 ,t^N),\\
\theta(x,\,y,\,z,\, t)\ =\ 1\hspace{0.2in}\mbox{on}\ \ {\Gamma}_{\ell}\times (0 , t^N),\\
 \theta(x,\,y,\,z,\, t)\ =\ 0\hspace{0.2in}\mbox{on}\ \ {\Gamma}_{r}\times (0 , t^N),\\
\frac{\D \partial \theta}{\D \partial \hat{ n}} ( x,\,y,\,z,\, t)\ =\ 0\hspace{0.2in}
\mbox{on}\ \ {\Gamma}\setminus \left( {\Gamma}_{l}\cup{\Gamma}_{r} \right) \times (0 , t^N).
\end{array}\right. \label{eq:NDicbc}
\end{equation}
In order to obtain the  dimensionless equations,  the distance
between the colder and hotter surfaces, $L_x$, has been chosen
as the reference length and the  scale factors for velocity, time
and  pressure are chosen as,  $\alpha / L_x$, ${L_x}^2 /\alpha$,
  $\rho_0 {\alpha}^2 /{L_x}^2$  respectively. The dimensionless
temperature is defined as $\theta\ =\ \frac{\D{T - T_c}}{\D{T_h
- T_c}}$.

Thus, the  three-dimensional rectangular model domain
has dimensions  $A_x\times A_y\times A_z$ where $A_{x_i}$ is the
\emph{aspect ratio}  in the $x_i$ direction.

\section{\normalsize Time Discretization by Marchuk-Yanenko-Type Operator Splitting Method of Problem  Equations(\ref{eq:NDNS}) -  (\ref{eq:NDicbc})}

 Let $\triangle t > 0 $ be  the time step and $t^n\ = \ n \triangle t$.
  At every time interval $[t^n,\, t^{n + 1}]$, the
Marchuk-Yanenko-type operator-splitting method involves a
sequence of computations as follows:
\begin{equation}\mbox{(I)}\hspace{0.2in} \mathbf{u}^0 =  \mathbf{u}_0\ ;\ \theta^0 = \theta_0,\hspace{3.5in}
\end{equation}
then, for $n \geq 0,$ $\mathbf{u}^n, \theta^n$ given, we compute
$\left( \theta ^{n + 1/3} , \mathbf{u}^{n + 1/3} , p^{n + 1}\right)$,
$\left( \theta ^{n + 2/3} , \mathbf{u}^{n + 2/3} \right)$  and  $
\left( \theta ^{n + 1} , \mathbf{u}^{n + 1} \right)$ as follows:

(II)\hspace{0.2in}Solve the \emph{pressure sub-problems}:\hspace{0.5in}
\begin{equation}\mbox{}\left \{
\begin{array}{l}
{\theta}^{n + 1/3} =  {\theta}^n \ \ \mbox{in}\ \Omega,\\
\frac{\D{ \mathbf{u}^{n + 1/3}\ -\ \mathbf{u}^n}}{\D {\triangle t}}\  =\ - \nabla p^{n + 1}\ \ \mbox{in}\ \Omega,\\
\mathbf{u}^{n + 1/3}\ =\ \mathbf{0} \ \ \mbox{on}\ \Gamma,\\
\nabla \cdot \ \mathbf{u}^{n + 1/3}\ =\ 0\ \ \mbox{in}\ \Omega, \label{eq:MYpres}
\end{array}\right.
\end{equation}
(III)\hspace{0.2in}then solve the \emph{transport sub-problems}:
\begin{equation}\left \{
\begin{array}{l}
\frac{\D\partial\theta}{\D\partial t} + {\mathbf{u}}^{n + 1/3}\cdot
\nabla {\theta}\ =\ 0 \ \ \mbox{in}\ \Omega\times (t^n ,\,
t^{n+1}),\\
\frac{\D \partial \mathbf{u}}{\D  \partial  t} + \mathbf{u}^{ n
+ 1/3}\cdot \nabla \mathbf{u}\ =\ \mathbf{0} \ \ \mbox{in}\
\Omega\times (t^n ,\,
t^{n+1}),\label{eq:MY_trans}\\
\mathbf{u}(t^n) = \mathbf{u}^{n + 1/3},\ \  \theta(t^n) =
\theta^{n + 1/3}.
\end{array}\right.
\end{equation}
\begin{equation}\mbox{Set} \hspace{0.2in} \mathbf{u}^{n + 2/3}  =
\mathbf{u}(t^{n + 1}), \hspace{0.2in}\theta^{n + 2/3} = \theta(t^{n+1}).\end{equation}
(IV)\hspace{0.2in}Finally, solve the \emph{diffusion   sub-problems}:
\begin{equation}\left \{
\begin{array}{l}
\frac{{\D {\theta}^{ n + 1} - {\theta}^{n +2/3}}}{\D {\triangle t}} -  \Delta {\theta}^{n + 1}\ =\ 0 \hspace{0.2in}\mbox{in}\ \Omega,\\
{\theta}^{n + 1}\ =\ 1\ \mbox{on}\ \ {\Gamma}_{\ell},\\
{\theta}^{n + 1}\ =\ 0\ \mbox{on}\ \ {\Gamma}_r,\\
\frac{\D{ \partial{\theta}} }{\D {\partial n}}^{n + 1}\ =\  0\ \ \mbox{on}\ \ {\Gamma}\setminus \left( {\Gamma}_l\cup {\Gamma}_r \right ),\\
\frac{\D {\mathbf{u}^{ n + 1} -\mathbf{u}^{ n + 2/3}}}{\D {\triangle t}} - Pr \Delta \mathbf{u}^{ n + 1}\ =\ Ra\, Pr\,{\theta}^{n + 1}\,\hat{\jmath} \ \ \mbox{in}\ \Omega,\\
\mathbf{u}^{ n + 1}\ =\ \mathbf{0} \ \ \mbox{on}\ \ \Gamma.\label{eq:MY_diff}
\end{array}\right.
\end{equation}

\section{\normalsize On  the Finite Element Approximation of sub-problems  (\ref{eq:MYpres}) -  (\ref{eq:MY_diff})}

The entire domain  $\Omega = \Omega_h$, constituting the computational domain,  is discretized into a
finite set, $\mathcal T_h$,  of tetrahedra  inside which
velocity, pressure and temperature are continuous and the
collection of the  tetrahedra satisfies the following properties:
\begin{enumerate}
\item $T_h = \overline{T_h} \subset {\overline{\Omega}}_h$ for all $T_h \,\in\,\mathcal{T}_h$.
\item $\mathcal{T}_h = \left\{ T_h \right\}$ is finite.
\item For $T_1\,,\, T_2 \in \mathcal{T}_h$,  if $\mbox{interior}(T_1)\  
\neq\  \mbox{interior}(T_2)$, then only one of the following is possible:
\begin{itemize}
\item $T_1 \cap T_2 = $ a vertex common to $T_1$
and $T_2$  or,
\item $T_1 \cap T_2 = $ a common edge or face  of  $T_1$ and $T_2$ or,
\item $T_1 \cap T_2 =  \emptyset$.
\end{itemize} \item $\bigcup\limits_{T_h \in \mathcal{T}_h} = {\overline{\Omega}}_h$.
\end{enumerate}

The  finite element  mesh for pressure is twice coarser than the
mesh for velocity. The mesh for temperature is the same as that
for velocity.  Let $\mathcal{T}_h$ be  the finite collection of
the tetrahedra for velocity and $\mathcal{T}_{2h}$, a similar
collection  for  pressure. Piecewise-linear approximation is
employed  for all variables including pressure. Hence,  the
vertices of the tetrahedra in $\mathcal{T}_h$ and
$\mathcal{T}_{2h}$ form  the \emph{nodes} for the finite element
meshes for velocity,temperature and pressure respectively.

A new  method of discretizing $\Omega_h$ and $\Omega_{2h}$ into
 4-node tetrahedral elements,  such that the pressure tetrahedra
are macro-elements consisting of eight sub-tetrahedra, is
 presented in this report. This is achieved by first
discretizing the domain into a finite set of $8$-cornered
\emph{brick-like} macros, each of which is then split into
two prisms along a vertical mid-plane through a diagonal
line on the top face and through the center of gravity of the
\emph{brick-like} macro. Each of the two prims is then divided
into $3$ tetrahedra in a very unique way.  $\mathcal{T}_{2h}$
is  constructed first  and,  as in two-dimensional cases,
$\mathcal{T}_{h}$ is constructed from $\mathcal{T}_{2h}$ by
connecting the edge midpoints on the faces and on any
resulting vertical mid-plane. A detailed explanation of the steps
involved in this \emph{tetrahedralization} is given in the
appendix. We however make the following remarks here:

\section{\normalsize Remarks}

\begin{enumerate}
\item The six generic tetrahedral elements have different shapes
but  the same volume which is one-sixth  of the volume of the
particular \emph{brick-like} macro containing them.
\item In order for each pressure tetrahedral  element to properly contain
exactly eight sub-tetrahedra, after the edge-midpoints have been
connected, it is necessary that, one set of three tetrahedra on
one prism must be a reflection of the second set of  three
tetrahedra on the  other prism about the vertical plane through
the diagonal of the \emph{brick-like} macro.
\item Although, a method, for splitting a  \emph{brick-like} macro into six
tetrahedra was given, by Zienkiewcz, in \cite{bib:zienkwcz77},
the discussion  did not include the situation of multi-level
grids as in this present case. \end{enumerate}

\section{\normalsize Discrete sub-problems and Weak Formulations}

Weak formulation of each set of sub-problems, determined by the
operator-splitting,  are obtained using the following
fundamental discrete spaces:
\begin{gather}
V_h =\left\{v_h \mid v_h \in C^0\left(\overline{\Omega}_h \right),\, {v_h}\vert_T \ \in\
\mathcal{P}_1,\, \forall\ T \in\  \mathcal{T}_h\right\},\label{eqn:femwave}\\
V_{0h} =\left\{v_h\mid v_h\, \in\, V_h,\, v_h\vert_{\Gamma} = 0 \right\},\label{eqn:femvel}\\
\Pi_{0h} =\left\{\phi_h\mid\phi \, \in\, V_h,\, \phi_h\vert_{\Gamma_\ell \cup \Gamma_r} = 0 \right\},\label{eqn:femtemp}\\
P_{h} = \left\{q_{h} \mid q_{h}\,\in\,C^0\left(\overline{\Omega}_{2h}\right )\,,\,q_{h}\vert_T\,\in\,\mathcal{P}_1\,,\
\forall T\,\in\,\mathcal{T}_{2h},\, \and \,\int_{\Omega_{2h}} q_h d\Omega = 0\right\}.\label{eqn:fempres}
\end{gather}
In equations(\ref{eqn:femwave}) - (\ref{eqn:fempres}),
$\mathcal{P}_1$ is the space of polynomials in three variables of
degree $\le 1$. The discrete approximation   associated to  the
finite element spaces described   above, for the weak formulation
of the pressure sub-problems  is: \begin{gather} \mbox{ }
\hspace{-0.5in}\mbox{Find  }  \{ \mathbf{u}_{h}^{n + 1/3},\,
p_{h}^{n + 1}\}\  \in\  \left(V_{0h}\right)^3\times P_{h} \mbox{
such that,}\nonumber\\ \frac{1}{\triangle t}\int_{\Omega_h}u^{n +
1/3}_i v_h \,d\Omega = \int_{\Omega_h}p_h^{n + 1}
\frac{\partial v_h}{\partial x_i} d\Omega +  \frac{1}{\triangle t}\int_{\Omega_h}u^n_i v_h \,d\Omega
\ \ \forall\ v_h\,\in\,V_{0h},\label{eqn:sqr1}\\
\mbox{  for each component  of }\mathbf{u}_h = \left\{u_i\right\}^3_{i = 1},\nonumber\\
\int_{\Omega_{2h}}q_h\,\nabla \cdot \mathbf{u}_h^{n +1/3} \,d \Omega\ =\
 0\, , \ \  \ \ \ \forall\ q_h\ \in\   P_h,\\
\mathbf{u}_h^{n + 1/3} = \mathbf{0},\hspace{0.2in}
\mbox{on}\ \  \Gamma.\label{eqn:soldiv} \end{gather}

The transport sub-problems combined with the wave equation
approach  described in \cite {bib:Mythesis} and
\cite{bib:DrPanDrGlo}   give a set of semi-discrete  sub-problems
 with weak formulation given as:
 
Find  $\theta_h\  \in\  V_{h} $ and  $\mathbf{u}_h\ \in\
 \left(V_{0h}\right)^3$,    $\forall t\,\in\, \left[t^n,\, t^{n +
1}\right]$ such that,
\begin{gather}
\mbox{} \hspace{-0.5in}\int_{\Omega_h}\left( \frac{\partial^2 \theta_h}{\partial t^2}\right)\phi_h\, d\Omega\ +
\int_{\Omega_h}\left(\mathbf{u}_h^{n + 1/3}\cdot\nabla \phi_h\right) \left(\mathbf{u}_h^{n + 1/3}
\cdot\nabla \theta_h\right)\, d\Omega = 0 ,\  \  \forall\ \phi_h\, \in\,    V_h, \label{eqn:wt6} \\
\mbox{} \hspace{-0.5in}\int_{\Omega_h}\left(\frac{ \partial^2\mathbf{u}_h}{\partial t^2}\right)\cdot\mathbf{w}_h\, d\Omega +
\int_{\Omega_h}\left( \mathbf{u}_h^{n + 1/3}\cdot\nabla \mathbf{w}_h\right) \left(\mathbf{u}_h^{n + 1/3} \cdot \nabla
\mathbf{u}_h\right)\, d \Omega = 0,\ \  \forall\ \mathbf{w}_h\,
\in\,  \left(V_h\right)^3,\label{eqn:wv6}\\
\theta_h (t^n ) =  \theta_h^{n + 1/3}, \hspace{0.8in}\mathbf{u}_h(t^n ) =   \mathbf{u}_h^{n + 1/3}, \\
\theta(t) = \left\{
\begin{array}{c}
1\hspace{0.2in}\mbox{on}\ \ {\Gamma}_{\ell},\\
0\hspace{0.2in}\mbox{on}\ \ {\Gamma}_{r}\label{eqn:Ineq1}
\end{array}\right.
\end{gather}

The solution of a  wave-like equation such as eq.(\ref{eqn:wt6})
and  eq.(\ref{eqn:wv6})  has been described in \cite
{bib:Mythesis} and  \cite{bib:DrPanDrGlo} for uniformly
structured meshes.  It involves time discretizing the equation by
 a second-order finite difference  scheme   with   an
initialization step consisting  of solution of  full discrete
version of eq.(\ref{eq:MY_trans}).   It is  noteworthy to mention
 that since  the discrete wave-like equations are explicit there
 is no need to store any square matrix thereby, conserving
computer memory. However, a local time   step $\triangle t/Q$ has
to be chosen with integer $Q$ sufficiently large so that the CFL
condition is not violated.

The set eq.(\ref{eq:MY_diff}),  of diffusion subproblems is
approximated by the following discrete sub-problems:
 
Find  $\theta^{n + 1}_h\  \in\ V_{h} $ and  $\mathbf{u}^{n +
1}_h\ \in\  \left( V_{0h}\right)^3$,  such that,
\begin{gather}
\int_{\Omega_h}\left( \frac{\theta^{n + 1}_h - \theta^{n +
2/3}_h}{\triangle t}\right)\phi_h\, d\Omega\ +
\int_{\Omega_h}\nabla \theta_h^{n + 1}\cdot\nabla \phi_h\,
d\Omega = 0 ,\  \  \forall\ \phi_h\hspace{0.2in}\in\    \Pi_{0h},
\label{eqn:diffA} \\ \theta_h^{n + 1}  =\left\{
\begin{array}{c}0\ \mbox{ on }\Gamma_r, \\  1\ \mbox{ on
}\Gamma_{\ell},\end{array} \right.
 \end{gather}
\begin{multline}
\int_{\Omega_h}\left( \frac{\mathbf{u}^{n + 1}_h - \mathbf{u}^{n
+ 2/3}_h}{\triangle t}\right)\cdot \mathbf{w}_h\, d\Omega\ +  Pr
\int_{\Omega_h}\nabla \mathbf{u}_h^{n + 1}\cdot\nabla
\mathbf{w}_h\, d\Omega \\ \mbox{} \hspace{0.8in}  =  Ra
Pr\int_{\Omega_h}\theta^{n + 1}_h\hat{\jmath} \cdot
\mathbf{w}_h\, d\Omega, \  \  \forall\ \mathbf{w}_h\
\hspace{0.2in}\in\  \left(V_h\right)^3,\label{eqn:diffB}\end{multline}

 
\section {\normalsize Solution Strategy for  the Pressure Subproblems (\ref{eqn:sqr1})  -  (\ref{eqn:soldiv})}
Let $\mathbf{x} = ( x,\,y,\,z)$ and $\D \left\{\omega_m
(\mathbf{x})\right\}_{m = 1}^{N_h}$  be the vector basis for
$V_{h}$ such that,
 \begin{equation}
  \omega_m (\mathbf{x}_{\ell}) = \left\{
 \begin{array}{c} 1\ \mbox{ if } \mathbf{x}_{\ell}  =  \mathbf{x}_m \\
 0\ \mbox{ if } \mathbf{x}_{\ell}  \ne \mathbf{x}_m
 \end{array}  \right.
 \end{equation}
 for  each  node on the  velocity mesh.  After applying Galerkin method  on the discrete weak formulation we get the following
 discrete sub-problem  at each node $\mathbf{x}_m$  on the
velocity mesh: \begin{gather} \mbox{ } \hspace{-0.5in}\mbox{Find
}\   \{ \mathbf{u}_{h}^{n + 1/3},\, p_{h}^{n + 1}\}\  \in\
\left(V_{0h}\right)^3\times P_{h} \mbox{  such that,}\nonumber\\
\frac{1}{\triangle t}\int_{\Omega_h}u^{n + 1/3}_i \omega_m
\,d\Omega =  \int_{\Omega_h}p_h^{n + 1}
\frac{\partial \omega_m}{\partial x_i} d\Omega + \frac{1}{\triangle t}\int_{\Omega_h}u^n_i \omega_m \,d\Omega,\   (1 \le i \le 3) \label{eqn:sqr2}\\
\mathbf{u}_h ^{n + 1/3}= \mathbf{0},\hspace{0.2in} \mbox{on}\
\Gamma,\\
\int_{\Omega_{2h}}\hspace{-2ex}\omega_k\,\nabla \cdot \mathbf{u}_h^{n +1/3} \,d \Omega\ =\  0,\, \mbox{ at each node }
\  \mathbf{x}_k   \mbox{   on the pressure mesh.}\label{eqn:soldivB}
\end{gather}
This sub-problems is solved by a preconditioned conjugate gradient (PCG) algorithm described in \cite{bib:Mythesis} and
\cite{bib:DrPanDrGlo}. Proper evaluation of  some integrals in
the PCG algorithm is  very crucial to   the overall performance of  this numerical technique. These include integrals of the forms  $\dint_{\Omega_h}p_h^{n + 1}\frac{\partial \omega_m}{\partial x_i} d\Omega$
and $\dint_{\Omega_{2h}}\hspace{-2ex}\omega_k\,\nabla \cdot \mathbf{u}_h^{n +1/3} \,d \Omega$
which involve product of  discrete functions over coarse and fine meshes. The following is a  summary of the techniques applied to these integrals:
 
\begin{enumerate}
\item All integrals of the form $\dint_{\Omega} fg\, d \Omega$  are approximated by the trapezoidal method globally on the element domain.
\item All integrals of the form  $\dint_{\Omega_h}p_h\frac{\partial \omega_m}{\partial x_i} d\Omega$ are
computed  on the fine velocity mesh element-by-element. On each pressure element, the  function
$p_h$ is interpolated linearly along element edges. Also, over
each velocity element, $p_h$ is approximated  by the average of
 its nodal values on the vertices.

\item The integrals of the form $\dint_{\Omega_{2h}}\hspace{-2ex}\omega_k\,\nabla \cdot
\mathbf{u}_h \,d \Omega$ are  computed over the coarse pressure
mesh element-by-element. In each pressure
tetrahedral element,  the nodal  values of the coarse-mesh  basis
 function $\omega_k$ on the vertices of
each of the $8$ included velocity elements are obtained by
linear interpolation along the  edges.  $\nabla \cdot
\mathbf{u}_h$ is piecewise-constant over each velocity element
since $\mathbf{u}_h$ is approximated by a piecewise-linear
function. \end{enumerate}

The PCG steps include the solution of a Neumann problem with
solution belonging to the space of functions with mean-value
zero. Since this problem is solved on a  relatively coarser
mesh,  with   banded storage,   it is expected that   memory
requirement will be a manageable size on most modern computers.
Thus  the  Neumann  problem    is solved by direct method after
cholesky factorization as suggested by Glowinski in \cite[page
267]{bib:glowinski91}.  However, since its solution has
mean-value zero the following  steps must be performed together
with the direct method:

\begin{enumerate}
\item Set one of the unknowns to zero and delete the
corresponding row and column; The N$\times$N linear system in
$\phi$ (say) will reduce to (N $-$ 1)$\times$(N $-$ 1) linear
system in $\phi^{\prime}$, \item Solve the reduced linear system
by direct method, for $\phi^{\prime}$, \item compute  mean value
of $\phi^{\prime} = \Dfrac{1}{meas(\Omega)}\dint \phi^{\prime}\,
d\Omega$, \item Set $\phi = \phi^{\prime} - (\mbox{ mean value of
}\phi^{\prime}$). ( $\int \phi\, d\Omega$ must be zero or in
practice$\int \phi\, d\Omega \leq \varepsilon $).
\end{enumerate}

\section{\normalsize Solution Strategy for the Diffusion Subproblems (\ref{eqn:diffA}) and (\ref{eqn:diffB})} 

The integrals in eq.(\ref{eqn:diffA}) and eq.(\ref{eqn:diffB})
are assembled over all the fine-mesh tetrahedral elements. With
proper ordering of nodes, the diffusion sub-problems for
temperature  and velocity components result in linear systems
of nodal values of the form
\begin{equation}  \mathbf{A}\mathbf{x} = \mathbf{f}\label{eqn:LinSys}\end{equation}
where  $\mathbf{A}$ is a symmetric positive definite, banded, sparse matrix.  The bandwidth of $\mathbf{A}$,
however grows very rapidly as the resolution is increased, so
 that memory requirement  becomes  prohibitively  large even
on supercomputers, despite  banded storage. Thus it is more
practicable to solve these linear systems  for temperature and
velocity components by iterative methods.  A careful observation
reveals that with piecewise-linear  approximations, $\mathbf{A}$
has only  $15$ non-zero diagonals in each case. Since  it
is also symmetric, a substantial amount of memory is freed
by storing only the main diagonal and the $7$ non-zero upper (or
lower) diagonals.
 
The associated Dirichlet boundary condition should be enforced in
a manner that preserves  the symmetry and sparse nature of
$\mathbf{A}$.  A method of achieving this is discussed by Stasa
in \cite[pp 59-61]{bib:stasa85}. The diffusion linear systems of
the form  eq.(\ref{eqn:LinSys})  are solved  by conjugate
gradient algorithm of Hestenes and Steifels (CGHS)  given for
example in \cite{bib:jenmckeown93}. This  algorithm involves
matrix-vector multiplications which are performed within the
bandwidth and in a manner that prevents ``\emph{fill-ins}''. Only
7 multiplication operations are required per row. To accelerate
and ensure convergence, preconditioning is usually necessary. A
way to achieve this is to ensure that Gerschgorin disks  for the
iterates  are concentric  by performing  \emph{symmetric
scaling},  where the main diagonal elements are scaled to unity
before commencing the iteration process.  The linear system for
temperature is solved first.  Its  converged value is applied to
solve the system for the velocity components. Further, since the
three linear systems for the segregated velocity components  are
independent, they are solved concurrently at each iteration step
of CGHS algorithm.  After scaling, the number  of iterations
required  for convergence of CGHS algorithm is usually a minute
fraction of the size of the  linear  system. For example, on a
$41^3$ velocity mesh the CGHS algorithm  for the linear system
for temperature converged in   $29$ to $41$ iterations depending
on the value of $Ra$. On the same mesh, the CGHS algorithm  for
the   linear system  for the segregated velocity components
converged in  $37$ iterations.

\section{\normalsize A numerical Example}

The numerical techniques presented in this report have   been applied to simulate  natural convection in
a cubical enclosure containing  air with prandtl number $0.71$ at $Ra$ in the range $10^3$ to $10^6$.  Initially, the enclosed fluid is stationary
 and the uniform temperature in the enclosure and its boundaries
is $T_c\mbox{ }^0C$.  Later, the surface  $\Gamma_{\ell}$ is
heated uniformly  to a  temperature $T_h \mbox{ }^0C$ while the
temperature on the surface $\Gamma_r$  is  held fixed at
$T_c\mbox{ }^0C$. These two surfaces are maintained at these
temperatures thereafter   while the remaining four surfaces of
the cube are  considered to  be perfectly thermally insulated.
The resulting density   variation within the confined fluid is
assumed to be small enough that the Boussinesq approximation is
valid.

\section{\normalsize Results and Discussion}

Steady state results were obtained for the numerical example  on
unstructured meshes with $41^3$ and $45^3$ velocity nodes. Steady
state solution is assumed when \begin{equation}
\Dfrac {\|u^{new}_i - u^{old}_i \|_2}{\|u^{new}_i\|_2}  \leq \varepsilon,
\end{equation}
where $u_i$ is one of the velocity components  and $\varepsilon$ is taken as $10^{-5}$.

Mesh of up to $41^3$ velocity nodes was used  for  $Ra = 10^3$
 and $10^4$ with  $\triangle t = 1/4000$. For $Ra = 10^5$ and $10^6$, when the boundary layer is
relatively thinner, non-uniform mesh of $45^3$ velocity nodes
 was used with $\triangle t = 1/9000$.  The time increment in the
re-discretization of the transport sub-problems was taken in each
 case as  $\tau = \frac{\triangle t}{10}$.  The construction on the
non-uniform mesh, where $13/16$ is taken as estimate for  the
boundary layer thickness, is summarized in Table \ref{tabl:D3NU}.

\begin{table}[hbt]\centering\footnotesize{
\begin{tabular}{|c|c|c|}\hline
$Ra = 10^5$ and $10^6$&sub-interval& N\b{o} of divisions (coarse mesh)\\ \cline{2-3}
$45^3$ fine mesh& $0 \leq x,\,y,\,z  \leq 3/16$&5\\ \cline{2-3}
&$3/16 \leq x,\,y,\,z  \leq 13/16$&12\\ \cline{2 - 3}
&$13/16 \leq x,\,y,\,z  \leq  1$&5\\ \hline
\end{tabular}
\caption{Construction of non-uniform  meshes  in  cubical
enclosure.}\label{tabl:D3NU} } \end{table}

Computations were performed on DEC Alpha PW500au, a single
processor, virtual memory machine and a  linux desk-top with
512MB core memory. Computations for each $Ra$ were started from
the initial conditions given in equation(\ref{eq:NDicbc}).  A
mesh with  $41^3$ velocity nodes  requires $70$MB of memory
while  a mesh  with $45^3$ velocity nodes  requires $99$MB of
memory.  An iteration in time, consisting of solution of the
pressure, transport and diffusion subproblems, takes average of
$1.64$ minutes   of CPU time. The PCG  algorithm  converged in
$7$ to $8$ iterations after initial transients.

 The \textit{Nusselt number}, a measure of the dimensionless heat transfer rate   across the isothermal
 walls, was computed on the hot wall ($x = 0$), in terms of the
\textit{overall Nusselt number},  $\overline{Nu}$ and
\textit{y-averaged Nusselt number}, $Nu_{av}(z)$, defined by  the
following equations:
\begin{equation}
\overline{Nu} =  - \int_{z = 0}^{z = A_z} \int_{ y = 0}^{ y = A_y}
\left. \frac{\partial \theta}{\partial x} \right|_{x = 0} dy\,dz .  \label{eqn:NUbar}
\end{equation}
\begin{equation}
{Nu}_{av} =   -  \int_{ y = 0}^{ y = A_y}
\left. \frac{\partial \theta}{\partial x} \right|_{x = 0} dy.   \label{eqn:NUav}
\end{equation}
For the cubical enclosure $A_y = 1.0$ and  $A_z = 1.0$
The local heat flux $\partial \theta / \partial x$, was approximated   by a second order
 forward-difference formula and the integrals were evaluated
using  the trapezoidal rule.
 
 Variations of all variables with respect to $z$ were investigated at each of the Rayleigh numbers
applied. These variations in the $z$-direction are however weaker
in magnitude than in other  directions.  In  Figure~\ref{fig:UymaxZ3456} and Figure~\ref{fig:NuZ3456} the
distributions of $\left(u_y\right)_{max}(z)$  at  $x = 0.5$
and the  y-averaged Nusselt number, $Nu_{av}(z)$ , are
illustrated at  each Rayleigh number.

\section{\normalsize Validation}

  A comparison of the  values of $\overline{Nu}$ and the
\textit{mean Nusselt number}  $Nu_{av}(0.5)$  with results of
Fusegi \emph{et al.\@}~\cite{bib:fusegi91} at each of the
Rayleigh numbers  used is given in Table \ref{tab1:Nutab}.
The results obtained using this present numerical procedure
are in good agreement with those of Fusegi \emph{et al.\@}.  Also, in
agreement with Fusegi \emph{et al.\@}~\cite{bib:fusegi91}, Janssen \emph{et
al.\@}~\cite{bib:janshen93} and Mallinson  \emph{et al.\@}~\cite{bib:mallinsondv77}, at  $Ra = 10^5$ and $10^6$, the
$z$-variations of $Nu_{av}(z)$  are apparent near the end walls
($z = 0$ and  $z = 1$), where it increases sharply.   Janssen \emph{et al.\@}~\cite{bib:janshen93} also reported a
value of $0.2585$ for $\overline{Nu} {Ra}^{- 1/4}$ at  $Ra =
10^6$. That is, $\overline{Nu}$ is $ 8.6396$ at this
Rayleigh number and this represents a difference of only $0.27\%$ of the value obtained using this
present numerical procedure (See Table \ref{tab1:NutabJ}).

\begin{table}[hbt]\centering  \footnotesize{
\begin{tabular}{|l|c|c|c|c|}\hline                                   
$Ra$&Quantity&Wave Equation & Fusegi \emph{et al.\@}\cite{bib:fusegi91} & \% Error  \\ \hline
$10^3$&$\overline{Nu}$&1.2466&1.085  &  12.96 \%  \\ \cline{2 - 5}
       &$Nu_{av}(0.5)$&1.2563&1.105  & 12.04\%\\ \hline
$10^4$&$\overline{Nu}$&1.9737 &2.10&  - 6.4\%  \\ \cline{2 - 5}
       &$Nu_{av}(0.5)$&2.1461& 2.302 &- 7.26\%\\ \hline
$10^5$& $\overline{Nu}$&4.2055  & 4.361& - 3.70\%  \\ \cline{2 - 5}
&$Nu_{av}(0.5)$&4.497&4.464& 0.73\% \\\hline
$10^6$& $\overline{Nu}$ &8.6628  & 8.770 &- 1.24\%  \\ \cline{2 - 5}
&$Nu_{av}(0.5)$&8.8434&9.012& - 1.91 \% \\ \hline
\end{tabular}
\caption{Comparison of Nusselt numbers with  results of Fusegi \emph{et al.\@} }\label{tab1:Nutab}
}
\end{table}
\begin{table}[hbt]\centering  \footnotesize{

\vspace{0.5in}

\begin{tabular}{|l|c|c|c|c|}\hline                                   
$Ra$&Quantity&Wave Equation & Janssen \emph{et al.\@}\cite{bib:janshen93} & \% Error  \\ \hline
 $10^6$& $\overline{Nu}{Ra}^{- 1/4}$ &   & $0.2585$&  \\
    & $\overline{Nu}$ &$8.6628$  & $ 8.6396$ &$0.27\%$ \\  \hline
\end{tabular}
\caption{Comparison of Nusselt numbers with  results of  Janssen \emph{et al.\@} }\label{tab1:NutabJ}
}
\end{table}

In agreement with  Janssen \emph{et
al.\@}~\cite{bib:janshen93}, at $Ra = 10^5$ and $10^6$,
 $\left(u_y\right)_{max}(z)$ has  two sharp peaks close to the
lateral walls. At $Ra = 10^5$, the peaks occur at $(0.075,\,
0.5,\, 0.15)$ and  $(0.075,\, 0.5,\, 0.85)$. At $ Ra = 10^6$, the peaks
occur at $(0.0375,\, 0.5,\, 0.0937)$ and   $(0.0375,\, 0.5,\,
0.9625)$. The value and location of the $z-$direction local maximum in each case are given in Table \ref{tab1:JansenC}.
\begin{table}[hbt]\centering  \footnotesize{

\vspace{0.5in}

\begin{tabular}{|l|c|c|}\hline
   & $z$-direction local maximum&Location at $y = 0.5$\\  \hline
$Ra = 10^5$&$89.442$ &$x = 0.075,\,z = 0.15$  \\ \hline
$ Ra = 10^6$&$314.1738$&$x = 0.0375,\,z = 0.0937$\\ \hline
\end{tabular}
\caption{Comparison of the peaks for $\left(u_y\right)_{max}(z)$ with  results of  Janssen \emph{et al.\@} }\label{tab1:JansenC}
}
\end{table}



The  transverse  variations  of characteristic quantities, were further investigated. From the mesh plots of  velocity
components, temperature and pressure on the $yz$-plane at $x =
0.5$ (not shown, see ref. \cite{bib:Mythesis}), it was observed that for a generic variable  $f$, representing
 $u_x$, $u_y$,  $\theta$ or $p$,  the following relation holds:
\begin{equation}f(x,\,y,\,z) \approx f(x,\,y,\,1 -  z)\  \  \forall\
x,\,y,\,z.\end{equation}
$u_z$ on the other hand satisfies the relation,
\begin{equation}
u_z(x,\,y,\,z) \approx - u_z(x,\,y,\,1 -  z)\  \  \forall\  x,\,y,\,z.
\end{equation}
This symmetry property is sometimes imposed by researchers in order to reduce calculations to only half of the entire
computational domain.

Also, from the contour plots on $xy$-plane (at $z = 0.5$) and
mesh plots  on $xz$-plane (at $y = 0.5$) (ref. \cite{bib:Mythesis}), it was observed that
there is also a special type of symmetry about the line
$(0.5,\,0.5,\,z)$ (that is, the line $y = 0.5$ on $yz$-plane at
$x = 0.5$), through the center of gravity, so that the following
relations are satisfied:
\begin{gather} u_x(x,\,y,\,z) \approx - u_x(1 -
x,\,1 - y,\,z),  \\ u_y(x,\,y,\,z) \approx - u_y(1 - x,\,1 -
y,\,z), \\ u_z(x,\,y,\,z) \approx  +  u_z(1 - x,\,1 - y,\,z),\\
p(x,\,y,\,z) \approx  +  p(1 - x,\,1 - y,\,z),\\
\theta(x,\,y,\,z) \approx  1 -\theta(1 - x,\,1 - y,\,z).
\end{gather} These two spatial symmetries have been observed by
Janssen \emph{et al.\@}~\cite{bib:janshen93} and they exploited
it by performing actual computations over only a quarter  of the
entire cubical enclosure.
 
  \newpage
 \newcommand{\envbib}[5]{#1, #2, \emph{#3}, \textbf{#4}, #5.}
 \newcommand{\dnvbib}[4]{#1, #2, \emph{#3}, #4.}
 \newcommand{\ijnmf}{Int. J. for Numerical Methods in Fluids}
  \newcommand{\ijhmt}{Int. J. Heat Mass Transfer}
  \newcommand{\cfd}{Computational Fluid Dynamics Review}

 \pagebreak
 
 \begin{center}\large\textbf{Appendix}\end{center}
 
 \appendix

\section{\normalsize A Systematic Method of Discretizing  A Three-Dimensional Domain into Four-Node Tetrahedral Elements}

Let $\sum_h$ and $\sum_{2h}$ be the set of all nodes on
$\Omega_h$ and $\Omega_{2h}$  respectively.  The first step in
the domain discretization is the construction of
$\mathcal{T}_{2h}$ which  is illustrated by  the following steps:

\begin{itemize}
\item Divide the three dimensional domain into \emph{brick-like macros},
\begin{equation}
\mbox{ }\hspace{-0.6in}\bigcup_{j = 1}^{N_{2h}}\left\{(x,\,y,\,z)| x_j \leq x \leq x_j +
\triangle x\,;\  y_j \leq  y \leq  y_j + \triangle y\,;\  z_j \leq z
\leq  z_j + \triangle z \right \},
\end{equation}
where $ N_{2h} = card\left(\sum_{2h}\right)$.

A typical  \emph{brick-like} macro, with corners labeled
A, B, C, D, E, F, G, H,  is shown in Figure~\ref{fig:d3macro}. Point A has coordinate $(x_j,\,y_j,\,z_j)$.

\item Divide each  \emph{brick-like}  macro into two prisms -
one with corners H, E, B, C, D, A and the other prism with
corners   H, E, F, B, C, G along diagonals $\left\{\overline{\mbox{HC}},\,
 \overline{\mbox{EB}}\right\}$ or $\left\{ \overline{\mbox{DG}},\, \overline{\mbox{AF}}\right\}$ (see
Figure~\ref{fig:pris3dA}).

\item Sub-divide the  first prism  into three tetrahedra, through
 line segments $\left\{ \overline{\mbox{DE}}\right.$, $\overline{\mbox{DB}}$,
 $\left. \overline{\mbox{HB}}\right\}$.\\
(\textit{top vertex-to-bottom edge construction})

\item Sub-divide the second prism  into three tetrahedra, through
 line segments $\left\{ \overline{\mbox{FH}}\right.$,  $\overline{\mbox{FC}}$,
 $\left. \overline{\mbox{HB}}\right\}$.\\
(\textit{bottom edge-to-top vertex construction})
\end{itemize}

Thus, each $8$-cornered \emph{brick-like} macro on the coarser
 mesh is divided into six tetrahedral elements for
pressure.

Observe that the three  tetrahedral elements on the first prism
may be obtained from those on  second prism   by reflecting
first prism about the line segment, $\overline{\mbox{HC}}$.
 Thus, we may reverse  the order of connecting points between
these two prisms. That is, we may sub-divide the  first prism
  into three tetrahedra   through line segments
$\left\{ \overline{\mbox{AH}}\right.$,  $\overline{\mbox{AC}}$,
$\left. \overline{\mbox{HB}}\right\}$ (\textit{bottom edge-to-top vertex
construction}).  The   second  prism  must  correspondingly be
sub-divided into three  tetrahedra through line segments
$\left\{ \overline{\mbox{GE}}\right.$, $\overline{\mbox{GB}}$,
 $\left. \overline{\mbox{HB}}\right\}$ (\textit{top vertex-to-bottom edge
construction}).  We may also use diagonal $\overline{\mbox{EC}}$
in place of $\overline{\mbox{HB}}$. If this reflection property
is taken into consideration, any appropriate combination  of line
segments, as described above,  will produce the identical set of
six tetrahedral elements.  This method of discretizing a
brick-like macro into six tetrahedra is unique in the sense that,
  the  six  tetrahedra are generic elements for  all the
tetrahedral elements in $\mathcal{T}_{2h}$ and $\mathcal{T}_h$.
The  six generic tetrahedra are shown in Figure~\ref{fig:d3fig2}.

\setlength{\unitlength}{0.5cm}
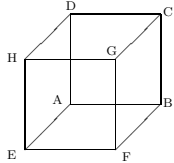
\begin{figure}[ht] \centering
\scalebox{0.6}{
\begin{picture}(4,6)(0,0)
\put(0,0){\line(1,0){4}}
\put(4,0){\line(0,1){4}}
\put(4,4){\line(-1,0){4}}
\put(0,4){\line(0,-1){4}}
\put(2,2){\line(1,0){4}}
\put(6,2){\line(0,1){4}}
\put(6,6){\line(-1,0){4}}
\put(2,6){\line(0,-1){4}}
\put(0,0){\line(1,1){2}}
\put(4,0){\line(1,1){2}}
\put(4,4){\line(1,1){2}}
\put(0,4){\line(1,1){2}}
\put(1.2,2){\FS A}
\put(6.1,1.9){\FS B}
\put(6.1,5.9){\FS C}
\put(1.8,6.2){\FS D}
\put(-0.8,-0.4){\FS E}
\put(4.3,-0.5){\FS F}
\put(3.6,4.2){\FS G}
\put(-0.8, 3.9){\FS H}
\end{picture}
}
\caption{A typical \emph{brick-like} macro.}\label{fig:d3macro}
\end{figure}

The finite set  $\mathcal{T}_h$ for velocity is next constructed
from   $\mathcal{T}_{2h}$  by connecting the edge-midpoints on
the faces of each tetrahedron and on any resulting vertical
mid-plane  following the steps for constructing the six generic
tetrahedra itemized above.

Thus, each pressure tetrahedral element is sub-divided into eight
 smaller tetrahedra, each of  which is similar in shape  to one
 of the six generic tetrahedra depicted in Figure~\ref{fig:d3fig2}.

Typical set  of $8$ velocity  elements in  two  generic pressure
macro-elements are delineated in Figure~\ref{fig:d3fig4}.

\setlength{\unitlength}{0.5cm}
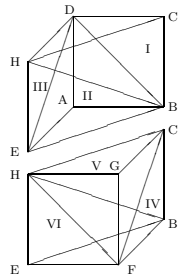
\begin{figure}[ht] \centering
\scalebox{0.6}{
\begin{picture}(4,11)(0,0)
\put(0,5){\line(1,1){2}}
\put(2,7){\line(1,0){4}}
\put(2,7){\line(0,1){4}}
\put(0,5){\line(0,1){4}}
\put(6,7){\line(0,1){4}}
\put(2,11){\line(1,0){4}}
\put(0,9){\line(1,1){2}}
\put(0,5){\line(3,1){6}} 
\put(0,9){\line(3,1){6}} 
\put(2,11){\line(1,-1){4}}
\put(0,5){\line(1,3){2}}
\put(0,9){\line(3,-1){6}}
\put(0,4){\line(3,-1){6}}
\put(0,0){\line(1,0){4}}
\put(4,0){\line(0,1){4}}
\put(4,4){\line(-1,0){4}}
\put(0,4){\line(0,-1){4}}
\put(6,2){\line(0,1){4}}
\put(4,0){\line(1,1){2}}
\put(4,4){\line(1,1){2}}
\put(0,4){\line(3,1){6}} 
\put(0,0){\line(3,1){6}} 
\put(0,4){\line(1,-1){4}}
\put(4,0){\line(1,3){2}}
\put(-0.8,-0.4){\FS E}
\put(-0.8,4.8){\FS E}
\put(4.4,-0.4){\FS F}
\put(6.2,6.8){\FS B}
\put(6.2,1.6){\FS B}
\put(-0.8,8.8){\FS H}
\put(-0.8,3.8){\FS H}
\put(3.6,4.2){\FS G}
\put(6.2,5.8){\FS C}
\put(6.2,10.8){\FS C}
\put(1.3,7.1){\FS A}
\put(1.6,11.1){\FS D}
\put(0.8,1.6){\FS {VI}}
\put(5.2,2.5){\FS  IV}
\put(2.8,4.2){\FS  V}
\put(0.2,7.7){\FS {III}}
\put(2.4,7.3){\FS {II}}
\put(5.2,9.4){\FS  I}
\end{picture}
}
\caption{A systematic division of a \emph{brick-like} macro into six tetrahedra.}\label{fig:pris3dA} 
\end{figure}

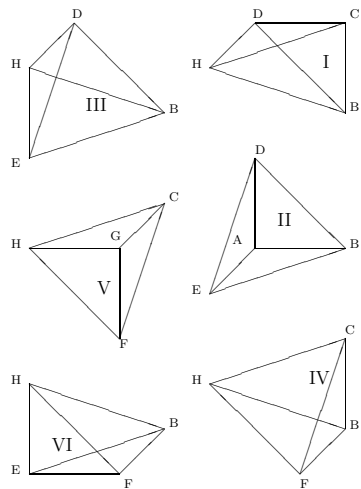
\begin{figure}[ht]\centering
\scalebox{0.6}{
\begin{picture}(14,22)(0,0)
\put(8,18){\line(3,1){6}} 
\put(8,18){\line(3,-1){6}}
\put(8,18){\line(1,1){2}}
\put(10,20){\line(1,0){4}}
\put(10,20){\line(1,-1){4}}
\put(14,16){\line(0,1){4}}
\put(13,18){{I}}
\put(14.2,16.1){\FS B}
\put(14.2,20.2){\FS C}
\put(9.9,20.2){\FS D}
\put(7.2,18){\FS H}
\put(10,14){\line(1,-1){4}}
\put(10,14){\line(0,-1){4}}
\put(10,10){\line(1,0){4}}
\put(8,8){\line(1,1){2}}
\put(8,8){\line(1,3){2}}
\put(8,8){\line(3,1){6}} 
\put(10,14.2){\FS D}
\put(9.0,10.2){\FS A}
\put(14.2,10){\FS B}
\put(7.2,8){\FS E}
  \put(11,11){{II}}
\put(0,18){\line(3,-1){6}} 
\put(0,18){\line(0,-1){4}}
\put(0,14){\line(1,3){2}}
\put(0,18){\line(1,1){2}}
\put(2,20){\line(1,-1){4}}
\put(0,14){\line(3,1){6}}
\put(-0.8,13.6){\FS E}
\put(6.2,16){\FS B}
\put(-0.8,18){\FS H}
\put(1.9,20.2){\FS D}
\put(2.5,16.1){{III}}
\put(14,2){\line(0,1){4}} 
\put(12,0){\line(1,1){2}}
\put(12,0){\line(1,3){2}}
\put(8,4){\line(3,1){6}}
\put(12,0){\line(-1,1){4}}
\put(8,4){\line(3,-1){6}}
\put(12,-0.6){\FS F}
\put(14,6.2){\FS C}
\put(14.2,2){\FS B}
\put(7.2,4){\FS H}
\put(12.4,4){{IV}}
\put(0,10){\line(3,1){6}} 
\put(0,10){\line(1,0){4}}
\put(4,6){\line(1,3){2}}
\put(4,6){\line(0,1){4}}
\put(4,10){\line(1,1){2}}
\put(4,6){\line(-1,1){4}}
\put(3,8){V}
\put(-0.8,10){\FS H}
\put(3.6,10.4){\FS G}
\put(4,5.6){\FS F}
\put(6.2,12.1){\FS C}
\put(0,4){\line(0,-1){4}}
\put(0,4){\line(1,-1){4}}
\put(0,4){\line(3,-1){6}} 
\put(0,0){\line(3,1){6}} 
\put(0,0){\line(1,0){4}}
\put(4,0){\line(1,1){2}}
\put(-0.8,0){\FS E}
\put(4.2,-0.6){\FS F}
\put(6.2,2.1){\FS B}
\put(-0.8,4){\FS H}
\put(1,1){{VI}}
\end{picture}
}
\caption{The six generic tetrahedral elements.}\label{fig:d3fig2}
\end{figure}

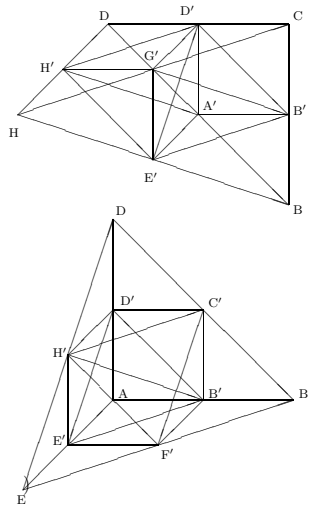
\begin{figure}[hbt]\centering

\scalebox{0.6}{
\begin{picture}(12,8)(0,0)
\put(0,4){\line(3,-1){12}}
\put(0,4){\line(3,1){12}}
\put(0,4){\line(1,1){4}}
\put(4,8){\line(1,0){8}}
\put(12,0){\line(0,1){8}}
\put(4,8){\line(1,-1){8}}
\put(2,6){\line(1,0){4}}
\put(6,6){\line(1,1){2}}
\put(2,6){\line(3,1){6}}
\put(6,2){\line(0,1){4}}
\put(6,2){\line(3,1){6}}
\put(6,6){\line(3,-1){6}}
\put(6,2){\line(1,1){2}}
\put(2,6){\line(1,-1){4}}
\put(2,6){\line(3,-1){6}}
\put(8,4){\line(1,0){4}}
\put(8,4){\line(0,1){4}}
\put(8,8){\line(1,-1){4}}
\put(6,2){\line(1,3){2}}
\put(12.2,4){\FS B$\mbox{}^{\prime}$}
\put(12.2,-0.4){\FS B}
\put(12.2,8.2){\FS C}
\put(5.6,1){\FS E$\mbox{}^{\prime}$}
\put(-0.4,3.0){\FS H}
\put(7.2,8.4){\FS D$\mbox{}^{\prime}$}
\put(8.2,4.2){\FS A$\mbox{}^{\prime}$}
\put(1.0 ,5.9){\FS H$\mbox{}^{\prime}$}
\put(5.6,6.4){\FS G$\mbox{}^{\prime}$}
\put(3.6,8.2){\FS D}
\end{picture}
}
\vspace{0.3in}

\scalebox{0.6}{
\begin{picture}(12,10)(0,0)
\put(0,0){\line(1,1){4}}
\put(4,4){\line(1,0){8}}
\put(4,4){\line(0,1){8}}
\put(0,0){\line(3,1){12}}
\put(0,0){\line(1,3){4}}
\put(4,12){\line(1,-1){8}}
\put(8,4){\line(-1,1){4}}
\put(4,8){\line(1,0){4}}
\put(8,4){\line(0,1){4}}
\put(4,8){\line(-1,-3){2}}
\put(4,8){\line(-1,-1){2}}
\put(2,2){\line(0,1){4}}
\put(2,2){\line(3,1){6}}
\put(2,6){\line(3,1){6}}
\put(8,4){\line(-1,-1){2}}
\put(2,2){\line(1,0){4}}
\put(2,6){\line(1,-1){4}}
\put(8,8){\line(-1,-3){2}}
\put(2,6){\line(3,-1){6}} )
\put(1.0,5.9){\FS H$\mbox{}^{\prime}$}
\put(1.0,2){\FS E$\mbox{}^{\prime}$}
\put(5.8,1.4){\FS F$\mbox{}^{\prime}$}
\put(7.9,4.1){\FS B$\mbox{}^{\prime}$}
\put(11.9,4.1){\FS B}
\put(3.9,4.1){\FS A}
\put(7.9,8.1){\FS C$\mbox{}^{\prime}$}
\put(4.0,8.2){\FS D$\mbox{}^{\prime}$}
\put(3.8,12.2){\FS D}
\put(-0.6,-0.6){\FS E}
\end{picture}
}
\caption{Typical pressure macro-elements  with $8$ sub-tetrahedral.}\label{fig:d3fig4}
\end{figure}

\begin{figure}[htb]\centering
\includegraphics[scale = 0.5]{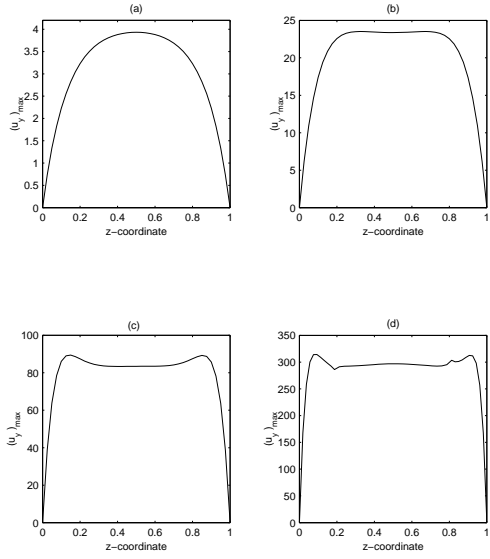}
\caption[Distribution of  $\left(u_y \right)_{max}$ along $z$-axis in cubical enclosure]{Distribution of $\left(u_y \right)_{max}$ along $z$-axis in  the cubical enclosure at\\ (a) $Ra = 10^3$  (b)  $Ra = 10^4$  (c)  $Ra = 10^5$  and (d) $Ra = 10^6$.}\label{fig:UymaxZ3456}
\end{figure}

\begin{figure}[htb]\centering
\includegraphics[scale = 0.5]{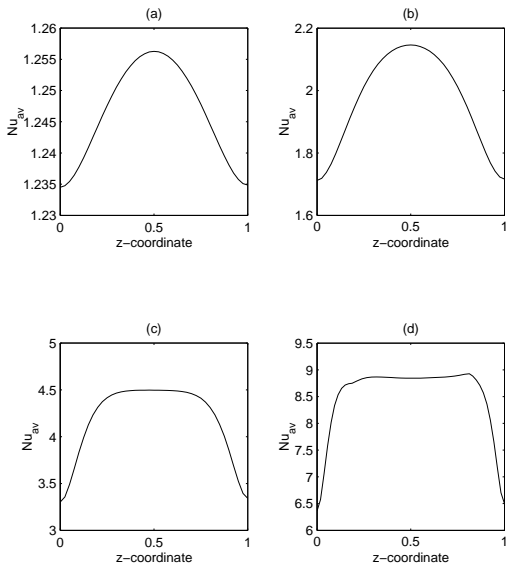}
\caption[Distribution of y-averaged Nusselt number,  $Nu_{av}$,
along $z$-axis  in the cubical  enclosure.]{Distribution of
y-averaged Nusselt number along $z$-axis in the cubical enclosure
at (a) $Ra = 10^3$  (b)  $Ra = 10^4$  (c)  $Ra = 10^5$  and (d)
$Ra = 10^6$.}\label{fig:NuZ3456} \end{figure}

\end{document}